\documentclass[12pt]{amsart}
\usepackage{amssymb}

\newtheorem{theorem}{Theorem}[section]
\newtheorem{proposition}[theorem]{Proposition}
\newtheorem{lemma}[theorem]{Lemma}
\newtheorem{corollary}[theorem]{Corollary}
\newtheorem{remark}[theorem]{Remark}
\newtheorem{definition}[theorem]{Definition}

\numberwithin{equation}{section}

\begin{document}

\title[Morse theory for Higgs $G$--bundles]{Morse theory 
for the space of Higgs $G$--bundles}

\author[I. Biswas]{Indranil Biswas}

\address{School of Mathematics, Tata Institute of Fundamental
Research, Homi Bhabha Road, Mumbai 400005, India}

\email{indranil@math.tifr.res.in}

\author[G. Wilkin]{Graeme Wilkin}

\address{Department of Mathematics, University of Colorado, Boulder, CO 80309-0395, USA}

\email{graeme.wilkin@colorado.edu}



\date{}

\begin{abstract}
Fix a $C^\infty$ principal $G$--bundle $E^0_G$ on a compact
connected Riemann 
surface $X$, where $G$ is a connected complex reductive linear
algebraic group. We consider the gradient flow of the
Yang--Mills--Higgs functional on the cotangent bundle of the space of 
all smooth connections on $E^0_G$. We prove that this flow preserves
the subset of Higgs $G$--bundles, and, furthermore, the flow
emanating from
any point of this subset has a limit. Given a Higgs $G$--bundle, we 
identify the limit 
point of the integral curve passing through it. These generalize the
results of the second named author on Higgs vector bundles.
\end{abstract}

\maketitle

\section{Introduction}\label{int}

The equivariant Morse theory of the Yang--Mills functional for vector 
bundles over a compact Riemann surface has been an extremely useful tool in 
studying the topology of the moduli space of semistable holomorphic 
bundles, beginning with the work of Atiyah and Bott in 
\cite{AtiyahBott83}, and continuing with the results of Kirwan on the 
intersection cohomology in \cite{Kirwan85-2}. In the case where the rank 
and degree of the bundle are coprime, this program was continued further 
by Kirwan in \cite{Kirwan92-2}, by Jeffrey and Kirwan (who computed the 
intersection pairings in \cite{JeffreyKirwan98}), and Earl and Kirwan in 
\cite{EarlKirwan04} who wrote down the relations in the cohomology ring.
As a result, via Morse theory we now have a complete description of the 
cohomology of this space.

The convergence properties of the gradient flow of the Yang--Mills 
functional were first studied by 
Daskalopoulos in \cite{Daskal92} and R{\aa}de in \cite{Rade92}. R{\aa}de
studies 
the more general problem of the Yang--Mills flow on the space of 
connections on a 2 or 3 dimensional manifold, and shows that the 
gradient flow converges to a critical point of the Yang--Mills 
functional. When the base manifold is a compact Riemann surface, then 
R{\aa}de's results show that there exists a Morse stratification of the 
space 
of holomorphic bundles. In \cite{Daskal92}, Daskalopoulos shows that 
this Morse stratification of the space of holomorphic bundles coincides 
with the algebraically defined Harder--Narasimhan stratification used 
by Atiyah and Bott, and uses this to obtain information about the 
homotopy type of the space of strictly stable rank $2$ bundles.

The analytically more complicated case of the Yang--Mills flow on the 
space of holomorphic structures on a K\"ahler surface (where bubbling 
occurs at isolated points on the surface in the limit of the flow) was 
studied by Daskalopoulos and Wentworth in \cite{DaskalWentworth04}, who 
showed that the algebraic and analytic stratifications coincide, and 
that the bubbling in the limit of the flow is determined by the 
algebraic properties of the initial conditions.

The main theorem of \cite{Wi} extend the gradient flow
results of Daskalopoulos and R{\aa}de to 
the space of Higgs vector bundles over a compact Riemann surface, where 
the functional in question is now the Yang--Mills--Higgs functional (see 
the definition below). These results were then used in \cite{DWWW} to 
carry out an analog of Atiyah and Bott's construction on the space of 
rank $2$ Higgs bundles, the first step in carrying out the 
Atiyah/Bott/Kirwan program described above for Higgs bundles. The 
purpose of the current paper is to generalize the gradient flow 
convergence theorem of \cite{Wi} to the space of Higgs principal bundles 
over a compact Riemann surface.

The main result of the paper can be stated as follows. Let $X$ be a 
compact Riemann surface. Given a $C^\infty$ vector bundle 
$V$ on $X$, the space of $C^\infty$ forms of type $(p\, ,q)$ with values
in $V$ will be denoted by $A^{p,q}(V)$. Let $G$ a connected reductive linear 
algebraic group over $\mathbb{C}$, and fix a maximal compact subgroup 
$K$. Fix a principal $K$--bundle $E_K^0$ with compact structure group 
$K$, and let $E_G^0$ denote the associated principal
bundle obtained by extending the 
structure group to $G$. Let $\mathcal{A}_0$ be the space of holomorphic 
structures on $E^0_G$, and consider the space $\mathcal{B}_0 = 
\mathcal{A}_0 \times A^{1,0}({\rm ad}(E_G^0))$ (we show in Section 
\ref{sec1} that this is the total space of the cotangent bundle of 
$\mathcal{A}_0$). A pair $(\overline{\partial}_{E_G^0}\, ,\theta)\,\in\, 
{\mathcal B}_0$ is called a \emph{Higgs pair} if $\theta$ is 
holomorphic with respect to $\overline{\partial}_{E_G^0}$. The space of 
Higgs pairs is denoted $\mathcal{S}(E_G^0)$, and the
\emph{Yang--Mills--Higgs functional} on the space $\mathcal{B}_0$ is 
given by
\begin{align*}
{\rm YMH}_G & : \mathcal{B}_0 \rightarrow \mathbb{R} \\
 {\rm YMH}_G (\overline{\partial}_{E_G^0} , \theta) & = \left\| 
K(\overline{\partial}_{E_G^0}) + [\theta, \theta^*] \right\|^2\, ,
\end{align*} 
where $K(\overline{\partial}_{E_G^0})$ is the curvature of the 
connection on 
$E_K^0$ associated to the holomorphic structure 
$\overline{\partial}_{E_G^0}$ 
on $E_G^0$ (see Section \ref{sec1} and the definition in equation 
\eqref{i} for full details of this construction).

The notion of the Harder--Narasimhan reduction of a Higgs principal 
bundle is recalled in Section \ref{sec:reduction}, and we show that the 
definition of socle reduction for semistable vector bundles (cf. 
\cite{HL}) extends to the semistable Higgs
$G$--bundles. Combining the Harder--Narasimhan reduction with the socle 
reduction gives the principal bundle analog of the graded object of the 
Harder--Narasimhan--Seshadri filtration studied in \cite{Wi}. The main 
theorem of this
paper generalizes the results of \cite{Wi} to Higgs principal bundles.

\begin{theorem} The gradient flow of ${\rm YMH}_G$ with initial 
conditions $(\overline{\partial}_{E_G^0}, \theta)$ in the space of 
Higgs pairs on 
$E_G^0$ converges to a Higgs pair isomorphic to the pair obtained by 
combining the socle reduction with the Harder--Narasimhan reduction of 
$(\overline{\partial}_{E_G^0}, \theta)$.
\end{theorem}

The idea of the proof is to reduce to the case of Higgs vector bundles 
studied in \cite{Wi}. Fix a $C^\infty$ principal $G$--bundle $E^0_G$
on a compact Riemann surface $X$.
Given a faithful representation $\rho \,:\, G 
\,\longrightarrow\, GL(V)$, let $\mathcal{S}(W)$ denote the space of 
Higgs pairs on the Hermitian vector bundle $W\, =\, E^0_G(V)$
associated to $E^0_G$ via $\rho$. We 
show that the Yang--Mills--Higgs flow on $\mathcal{S}(E_G^0)$ is induced 
by the Yang--Mills--Higgs flow on $\mathcal{S}(W)$, and that the 
convergence of the flow on $\mathcal{S}(W)$ (guaranteed by the results 
of \cite{Wi}) implies that the flow of ${\rm YMH}_G$ on 
$\mathcal{S}(E_G^0)$ converges also. The result then follows by showing 
that the Harder--Narasimhan--socle reduction on $W$ induces the 
Harder--Narasimhan--socle reduction on $E_G^0$.

The paper is organized as follows. Section \ref{sec1} contains the basic 
results necessary to the rest of the paper: the construction of the map 
$\phi : \mathcal{S}(E_G^0) \hookrightarrow \mathcal{S}(W)$ in equation 
\eqref{e9}, and a proof that the Yang--Mills--Higgs flows coincide
(Corollary \ref{cor1}). The main result of Section \ref{sec:closed} is 
Theorem \ref{lw}, which shows that the flow of ${\rm YMH}_G$ on 
$\mathcal{S}(E_G^0)$ converges. The results of the final section relate 
the Harder--Narasimhan--socle reduction of a Higgs pair 
$(\overline{\partial}_{E_G^0}\, , \theta) \in \mathcal{S}(E_G^0)$ to 
that of the associated Higgs pair
$\phi(\overline{\partial}_{E_G^0}\, ,\theta)\,\in\, \mathcal{S}(W)$.

\section{The Yang--Mills--Higgs functional}\label{sec1}

Let $G$ be a connected reductive linear algebraic group defined
over $\mathbb C$. Fix a maximal compact
subgroup
\begin{equation}\label{e1}
K\, \subset\, G\, .
\end{equation}
Fix a faithful representation
\begin{equation}\label{e2}
\rho\, :\, G\, \longrightarrow\, \text{GL}(V)\, ,
\end{equation}
where $V$ is a finite dimensional complex vector space. Fix a 
maximal compact subgroup $\widetilde{K}$
\begin{equation}\label{e3}
\rho(K) \, \subset\,\widetilde{K}\, \subset\, \text{GL}(V)
\end{equation}
of $\text{GL}(V)$.

The Lie algebra of $G$ will be denoted by
$\mathfrak g$. The group $G$ has the adjoint
action on $\mathfrak g$. So $\mathfrak g$ is
a $G$--module.

Let $X$ be a compact connected Riemann surface.
Fix a $C^\infty$ principal $K$--bundle
\begin{equation}\label{e4a}
E^0_K\, \longrightarrow\, X\, .
\end{equation}
Let
\begin{equation}\label{e4}
E^0_G\, :=\, E^0_K(G)\,=\, E^0_K\times_K G
\, \longrightarrow\, X
\end{equation}
be the principal $G$--bundle obtained by extending the structure
group of $E^0_K$ using the inclusion map $K\, \hookrightarrow\, G$.
Let $\text{ad}(E^0_G)\, :=\, E^0_G({\mathfrak g})\,=\,
E^0_G\times_G {\mathfrak g}$ be the adjoint bundle of $E^0_G$.
In other words, $\text{ad}(E^0_G)$ is the vector bundle over $X$
associated to the principal $G$--bundle $E^0_G$ for
the $G$--module $\mathfrak g$.

Let
\begin{equation}\label{e5}
{\mathcal A}_0\, :=\, {\mathcal A}_{E^0_G}
\end{equation}
be the space of all holomorphic structures on the principal
$G$--bundle $E^0_G$. We note that ${\mathcal A}_0$ is an affine
space for the vector space $A^{0,1}(\text{ad}(E^0_G))$, which is
the space of all smooth $(0\, ,1)$--forms with values in
$\text{ad}(E^0_G)$. Fix a holomorphic structure 
\begin{equation}\label{e6}
\overline{\partial}_0\, :=\, \overline{\partial}^0_{E^0_G}
\end{equation}
on $E^0_G$. Using $\overline{\partial}_0$, the affine space
${\mathcal A}_0$ gets identified with $A^{0,1}(\text{ad}(E^0_G))$.

Let $E_G\, \longrightarrow\, X$ be a holomorphic principal $G$--bundle.
A \textit{Higgs field} on $E_G$ is a holomorphic section of
$\text{ad}(E_G)\otimes K_X$ over $X$. A pair $(E_G\, ,\theta)$, where
$\theta$ is a Higgs field on $E_G$, is called a
\textit{Higgs $G$--bundle}.
A holomorphic structure on a principal $G$--bundle
$E_G$ defines a holomorphic structure
on the vector bundle $\text{ad}(E_G)\otimes K_X$. The Dolbeault
operator on $\text{ad}(E^0_G)\otimes K_X$ corresponding to any
$\overline{\partial}_{E_G^0}\, \in\, {\mathcal A}_0$ (see \eqref{e5})
will also be denoted by $\overline{\partial}_{E_G^0}$.
We note that a pair $(\overline{\partial}_{E_G^0}\, ,\theta)\, \in\,
{\mathcal A}_0\times A^{1,0}(\text{ad}(E^0_G))$ with
$\overline{\partial}_{E_G^0}(\theta)\,=\, 0$ defines a Higgs 
$G$--bundle.

Define
\begin{equation}\label{e7}
{\mathcal B}_0\, :=\,
{\mathcal A}_0\times A^{1,0}(\text{ad}(E^0_G))\, .
\end{equation}
So if $\overline{\partial}_{E_G^0}\, \in\, {\mathcal A}_0$,
and $\theta$ is a Higgs field on the holomorphic principal
$G$--bundle $(E^0_G\, ,\overline{\partial}_{E_G^0})$, then
$(\overline{\partial}_{E_G^0}\, ,\theta)\, \in\, {\mathcal B}_0$.
We will see later that ${\mathcal B}_0$ is the
total space of the
cotangent bundle of the affine space ${\mathcal A}_0$.

Let
\begin{equation}\label{W}
W\, :=\, E^0_K(V)\,\longrightarrow\, X
\end{equation}
be the vector bundle associated to the
principal $K$--bundle $E^0_K$ (see \eqref{e4a}) for the $K$--module
$V$ in \eqref{e2}. Therefore, $W$ is identified with the
vector bundle associated to the
principal $G$--bundle $E^0_G$ in \eqref{e4} for the $G$--module
$V$. A holomorphic structure on $E^0_G$ defines a holomorphic structure
on the vector bundle $W$. Using the injective homomorphism of Lie
algebras associated to $\rho$ in \eqref{e2}
\begin{equation}\label{l1}
{\mathfrak g}\, \longrightarrow\, \text{End}_{\mathbb C}(V)\, ,
\end{equation}
we get a homomorphism of vector bundles
\begin{equation}\label{h0}
\text{ad}(E^0_G)\, \longrightarrow\, End(W)\, =\, W\otimes W^*\, ,
\end{equation}
where $W$ is the vector bundle in \eqref{W}. Take any
$(\overline{\partial}_{E_G^0}\, ,\theta)\,\in\, {\mathcal B}_0$ (see
\eqref{e7}). Let $\overline{\partial}_W$ be the holomorphic
structure on $W$ defined by $\overline{\partial}_{E_G^0}$. Let
$\theta_W\, \in\, A^{1,0}(End(W))$ be the smooth section given by
$\theta$ using the homomorphism in \eqref{h0}.

Let ${\mathcal A}(W)$ be the space of all holomorphic
structures on the vector bundle $W$. Define
\begin{equation}\label{e8}
{\mathcal B}_W\, :=\,{\mathcal A}(W)\times A^{1,0}(End(W))\, .
\end{equation}
Since $W$ is associated to $E^0_G$ by a faithful
representation, there is a natural embedding
\begin{equation}\label{eb1}
\delta\, :\, {\mathcal A}_0\, \longrightarrow\, {\mathcal A}(W)\, ,
\end{equation}
where ${\mathcal A}_0$ is defined in \eqref{e5}.
We have an embedding
\begin{equation}\label{e9}
\phi\, :\, {\mathcal B}_0\, \longrightarrow\, {\mathcal B}_W
\end{equation}
that sends any $(\overline{\partial}_{E_G^0}\, ,\theta)$ to the
pair $(\overline{\partial}_W\, ,\theta_W)$ constructed
above from $(\overline{\partial}_{E_G^0}\, ,\theta)$.

The Lie algebra of $\widetilde K$ (see \eqref{e3}) will be
denoted by $\widetilde{\mathfrak k}$. Let $g_0$ denote the
inner product on $\widetilde{\mathfrak k}$ defined by
$\langle A\, ,B\rangle\, =\, -\text{trace}(AB)$.
Since $\text{End}_{\mathbb C}(V)\,=\, \widetilde{\mathfrak k}
\oplus \sqrt{-1}\widetilde{\mathfrak k}$, where $V$ is the
$G$--module in \eqref{e2}, this $g_0$ defines
a Hermitian inner product $g_1$ on $\text{End}_{\mathbb C}(V)$.
The Lie algebra of $K$ (see \eqref{e1}) will be
denoted by ${\mathfrak k}$. Let $g'_0$ be the restriction of
$g_0$ to the subspace $\mathfrak k\, \subset\,
\widetilde{\mathfrak k}$. Since $\mathfrak g\,=\,
{\mathfrak k}\oplus \sqrt{-1}{\mathfrak k}$, this $g'_0$ defines
a Hermitian inner product $g'_1$ on $\mathfrak g$. Note that
$g'_1$ is the restriction of $g_1$.
The inner products $g'_1$ and $g_1$ induce
inner products on the fibers of the vector bundle
$\text{ad}(E^0_G)$ and $End(W)$ respectively. Indeed, these
follow from the fact that $g'_1$ and $g_1$ are $K$--invariant
and $\widetilde K$--invariant respectively.

We will now show that the Cartesian product
${\mathcal B}_0$ in \eqref{e7} is the total space of
the cotangent bundle of the affine space ${\mathcal A}_0$.
For any $(\omega_{0,1}\, ,\omega_{1,0})\, \in\,
A^{0,1}(\text{ad}(E^0_G))\times A^{1,0}(\text{ad}(E^0_G))$,
we have
$$
\langle\omega_{0,1}\, ,\omega_{1,0}\rangle\, \in\,
A^{1,1}
$$
using the inner product on the fibers of $\text{ad}(E^0_G)$.
Consider the pairing
$$
A^{0,1}(\text{ad}(E^0_G))\times A^{1,0}(\text{ad}(E^0_G))
\,\longrightarrow\, \mathbb C
$$
defined by
$$
(\omega_{0,1}\, ,\omega_{1,0})\, \longmapsto\,
\int_X \langle\omega_{0,1}\, ,\omega_{1,0}\rangle\, \in\,
\mathbb C\, .
$$
This pairing identifies $A^{1,0}(\text{ad}(E^0_G))$ with the dual
of $A^{0,1}(\text{ad}(E^0_G))$. Therefore,
the total space of the cotangent bundle of the affine space
${\mathcal A}_0$ gets identified with ${\mathcal B}_0$.

Similarly, using the inner product on the fibers of $End(W)$, the
Cartesian product ${\mathcal B}_W$ in \eqref{e8} gets identified with 
the total space of the cotangent bundle of the affine space ${\mathcal 
A}(W)$.

Since the fibers of the vector bundles
$\text{ad}(E^0_G)$ and $End(W)$ have inner products,
we get inner products on the vector spaces
$$
A^{0,1}(\text{ad}(E^0_G))\oplus
A^{1,0}(\text{ad}(E^0_G)) ~\,~\, \text{~and~}~\,~\,
A^{0,1}(End(W))\oplus A^{1,0}(End(W))\, .
$$
More precisely, the inner product on $A^{0,1}(\text{ad}(E^0_G))\oplus
A^{1,0}(\text{ad}(E^0_G))$ is defined by
$$
\Vert(\omega_{0,1}\, ,\omega_{1,0})\Vert^2 \, =\, \sqrt{-1} \int_X
\langle \omega_{1,0}\, , \overline{\omega_{1,0}}\rangle -\sqrt{-1}
\int_X \langle \omega_{0,1}\, , \overline{\omega_{0,1}}\rangle\, .
$$
The inner product on $A^{0,1}(End(W))\oplus A^{1,0}(End(W))$ is
defined similarly.

Recall that ${\mathcal B}_0$ and ${\mathcal B}_W$ are identified
with
$$
A^{0,1}(\text{ad}(E^0_G))\oplus
A^{1,0}(\text{ad}(E^0_G)) ~\,~\, \text{~and~}~\,~\, 
A^{0,1}(End(W))\oplus A^{1,0}(End(W))
$$
respectively (the affine space ${\mathcal A}_0$
is identified with $A^{0,1}(\text{ad}(E^0_G))$ after choosing
the base point $\overline{\partial}_0$ in \eqref{e6}; since
$\overline{\partial}_0$ gives a point in ${\mathcal A}(W)$, it
follows that ${\mathcal A}(W)$ is identified with
$A^{0,1}(End(W))$). Therefore, the
inner products on $A^{0,1}(\text{ad}(E^0_G))\oplus A^{1,0}
(\text{ad}(E^0_G))$ and $A^{0,1}(End(W))\oplus A^{1,0}(End(W))$
define K\"ahler structures on ${\mathcal B}_0$ and ${\mathcal B}_W$
respectively.

\begin{lemma}\label{lem1}
The embedding $\phi$ in \eqref{e9} preserves the
K\"ahler forms. Moreover, the second fundamental form
of the embedding $\phi$ vanishes. In particular, this
embedding is totally geodesic.
\end{lemma}

\begin{proof}
Since the inner product $g'_1$ on $\mathfrak g$ is the restriction
of the inner product $g_1$ on $\text{End}_{\mathbb C}(V)$, it follows 
immediately
that $\phi$ preserves the K\"ahler forms.

Let ${\mathfrak g}^\perp\, \subset\, \text{End}_{\mathbb C}(V)$ be the
orthogonal complement for the inner product $g_1$ on
$\text{End}_{\mathbb C}(V)$. Since $g_1$ is $K$--invariant (recall that
it is in fact $\widetilde K$--invariant), and the adjoint action of
$K$ on $\text{End}_{\mathbb C}(V)$ preserves the subspace ${\mathfrak 
g}\, 
\subset\, \text{End}_{\mathbb C}(V)$, it follows that the adjoint action 
of
$K$ on $\text{End}_{\mathbb C}(V)$ preserves ${\mathfrak g}^\perp$.
Since $K$ is Zariski dense in $G$, it follows that the
adjoint action of $G$ on $\text{End}_{\mathbb C}(V)$ preserves 
${\mathfrak 
g}^\perp$. Therefore, the orthogonal decomposition
\begin{equation}\label{od2}
\text{End}_{\mathbb C}(V)\,=\, {\mathfrak g}\oplus {\mathfrak g}^\perp
\end{equation}
is preserved by the adjoint action of $G$.

Let
\begin{equation}\label{f0}
F_0\, :=\, E^0_K({\mathfrak g}^\perp)\, \longrightarrow\, X
\end{equation}
be the vector bundle associated to the principal $K$--bundle
$E^0_K$ (see \eqref{e4a}) for the $K$--module ${\mathfrak g}^\perp$.
Note that the $G$--invariant orthogonal decomposition of
$\text{End}_{\mathbb C}(V)$ in \eqref{od2} induces an orthogonal 
decomposition
\begin{equation}\label{od}
End(W)\, =\, \text{ad}(E^0_G)\oplus F_0\, .
\end{equation}
Hence we have orthogonal decompositions
\begin{equation}\label{o1}
A^{0,1}(End(W)) = A^{0,1}(\text{ad}(E^0_G))\oplus
A^{0,1}(F_0)\, \text{~and~}\,
A^{1,0}(End(W)) = A^{1,0}(\text{ad}(E^0_G))\oplus
A^{1,0}(F_0)\, .
\end{equation}

Let $\mathcal H$ denote the trivial vector bundle over ${\mathcal
B}_0$ (see \eqref{e7}) with fiber $A^{0,1}(F_0)\oplus A^{1,0}(F_0)$.
Using the orthogonal decompositions in \eqref{o1}
it follows that the orthogonal complement of the differential
\begin{equation}\label{e10}
d\phi\, :\, T^{1,0}{\mathcal B}_0\, \longrightarrow\,
\phi^*T^{1,0}{\mathcal B}_W
\end{equation}
is identified with the above defined vector bundle $\mathcal H$
(here $T^{1,0}$ denotes the holomorphic tangent bundle).

On the other hand, $\mathcal H\, \subset\, \phi^*T^{1,0}
{\mathcal B}_W$ is a holomorphic subbundle because the
adjoint action of $G$ on $\text{End}_{\mathbb C}(V)$ preserves 
${\mathfrak g}^\perp$. Consequently, the orthogonal
complement ${\mathcal H}\, =\, d\phi(T^{1,0}{\mathcal B}_0)^\perp
\, \subset\, \phi^*T^{1,0}{\mathcal B}_W$ (see \eqref{e10})
is preserved by the Chern connection on the holomorphic
Hermitian vector bundle $\phi^*T^{1,0}{\mathcal B}_W$. Since
the Chern connection for a K\"ahler metric coincides with
the Levi--Civita connection, it follows that
$d\phi(T^{1,0}{\mathcal B}_0)^\perp$ is preserved by the
connection on $\phi^*T^{1,0}{\mathcal B}_W$ obtained by
pulling back the Levi--Civita connection on
${\mathcal B}_W$. In other words, the
second fundamental form of the embedding $\phi$ vanishes.
This completes the proof of the lemma.
\end{proof}

\begin{remark}\label{rem1}
{\rm Take any $z\, :=\,
(\overline{\partial}_{E_G^0}\, ,\theta)\,\in\, {\mathcal B}_0$
which is a Higgs $G$--bundle, meaning
$\theta$ is holomorphic with respect to the holomorphic
structure $\overline{\partial}_{E_G^0}$. The $\phi(z)$ is a
Higgs vector bundle.}
\end{remark}

A connection on the principal $G$--bundle
$E^0_G$ decomposes the real tangent bundle
$T^{\mathbb R}E^0_G$ into a direct sum of horizontal and vertical
tangent bundles. Using this decomposition, the almost
complex structures of $G$ and $X$ together produce an almost
complex structure on $E^0_G$. Let $\overline{\partial}_{E_G^0}$ be
a holomorphic structure on $E^0_G$. A connection $\nabla$ on $E^0_G$
is said to be \textit{compatible} with $\overline{\partial}_{E_G^0}$
if the almost complex structure on $E^0_G$ given by $\nabla$
coincides with the one underlying the complex structure
$\overline{\partial}_{E_G^0}$ on $E^0_G$.

Given a holomorphic structure $\overline{\partial}_{E_G^0}$ on
$E^0_G$, there is a unique connection $\nabla$ on $E^0_K$
such that the connection on $E^0_G$ induced by $\nabla$
is compatible with $\overline{\partial}_{E_G^0}$; it is known
as the \textit{Chern connection}. On the other
hand, given a connection $\nabla_1$ on $E^0_G$, there is a
unique holomorphic structure $\overline{\partial}'_1$ on $E^0_G$
such that $\nabla_1$ is compatible with $\overline{\partial}'_1$
(this is because $\dim_{\mathbb C}X\,=\, 1$).
Therefore, we have a canonical bijective correspondence between 
${\mathcal A}_0$ (see \eqref{e5}) and the space of all connections
on $E^0_K$. Similarly, we have a canonical bijective correspondence 
between ${\mathcal A}(W)$ (the space of all holomorphic
structures on the vector bundle $W$ in \eqref{W}) and
the space of all connections
on the principal $\widetilde{K}$--bundle
\begin{equation}\label{ekt}
E^0_{\widetilde{K}}\, :=\, E^0_K(\widetilde{K})\,=\,
E^0_K\times_K \widetilde{K}\,\longrightarrow\, X
\end{equation}
obtained by extending the structure group of $E^0_K$ using
the homomorphism $\rho\, :\, K\, \longrightarrow\,
\widetilde{K}$.

The curvature of a connection $\nabla$ will be denoted
by $K(\nabla)$.

Let
\begin{equation}\label{st}
*\, :\, End(W)\, \longrightarrow\, End(W)
\end{equation}
be the
conjugate linear automorphism that acts on the subbundle
$\text{ad}(E^0_{\widetilde{K}}) \, \subset\, End(W)$
(see \eqref{ekt}) as multiplication by $-1$; since
$End(W)\,=\, \text{ad}(E^0_{\widetilde{K}})\oplus
\sqrt{-1}\cdot \text{ad}(E^0_{\widetilde{K}})$, this condition
uniquely determines the automorphism in \eqref{st}.

Fix a Hermitian metric $h_0$ on $T^{1,0}X$.

Let
$$
\text{YMH}_W\, :\, {\mathcal B}_W\, \longrightarrow\,
{\mathbb R}
$$
be the function defined by $(\overline{\partial}_{E_G^0}\, ,\theta)
\,\longmapsto\, \Vert K(\overline{\partial}_{E_G^0})+[\theta\, ,
\theta^*]\Vert^2$, where $K(\overline{\partial}_{E_G^0})$ is the
curvature of the connection associated to
$\overline{\partial}_{E_G^0}$, and the inner product on $2$--forms
is defined using $h$ and the inner product on the
fibers of $End(E)$. If locally $\theta\, =\, A\times dz$, then
$[\theta\, , \theta^*]\,=\, (AA^*-A^*A)dz{\wedge}d\overline{z}$;
see \cite{Wi} for more on this function $\text{YMH}_W$.

Let $*\, :\, \text{ad}(E^0_G)\, \longrightarrow\, \text{ad}(E^0_G)$
be the conjugate linear automorphism that acts on the subbundle
$\text{ad}(E^0_K) \, \subset\, \text{ad}(E^0_G)$ (see \eqref{e4a})
as multiplication by $-1$. Note that this automorphism coincides
with the restriction of the automorphism in \eqref{st}.
Consider ${\mathcal B}_0$ defined in \eqref{e7}. Let
\begin{equation}\label{e11}
\text{YMH}_G\, :\, {\mathcal B}_0\, \longrightarrow\,
{\mathbb R}
\end{equation}
be the function defined by $(\overline{\partial}_{E_G^0}\, ,\theta)
\,\longmapsto\, \Vert K(\overline{\partial}_{E_G^0})+[\theta\, ,
\theta^*]\Vert^2$; as before, $K(\overline{\partial}_{E_G^0})$ is the
curvature of the connection associated to the
holomorphic structure $\overline{\partial}_{E_G^0}$.

We first note that
\begin{equation}\label{i}
\text{YMH}_G\,=\, \text{YMH}_W\circ\phi\, ,
\end{equation}
where $\phi$ is the function constructed in \eqref{e9}. Let
$d\text{YMH}_W$ be the smooth exact $1$--form on ${\mathcal B}_W$.
The following lemma shows that the normal vectors to
$T^{\mathbb R}{\mathcal B}_0$ for the embedding $\phi$ in
\eqref{e9} are annihilated by the form $d\text{YMH}_W$.

\begin{lemma}\label{lem2}
For any point $x\, \in\, {\mathcal B}_0$, and any normal vector
$$
v\, \in\, (d\phi(T^{\mathbb R}_x{\mathcal B}_0))^\perp
\, \subset\, T^{\mathbb R}_{\phi(x)}{\mathcal B}_W\, ,
$$
the following holds:
$$
d{\rm YMH}_W(v)\, =\, 0\, .
$$
\end{lemma}

\begin{proof}
Take any pair $(v\, ,w)\, \in\, A^{0,1}(F_0)\oplus A^{1,0}(F_0)$,
where $F_0$ is the vector bundle in \eqref{f0}.
Take any $(\overline{\partial}_{E_G^0}\, ,\theta)\,\in\, {\mathcal 
B}_0$. Let $\nabla$ be the connection on $E^0_G$ corresponding to the
holomorphic structure $\overline{\partial}_{E_G^0}$. Therefore,
the connection on $E^0_G$ corresponding to the
holomorphic structure $\overline{\partial}_{E_G^0}+tv$, where $t\,\in\,
\mathbb R$, is $\nabla+tv - tv^*$. The automorphism in \eqref{st}
preserves the orthogonal decomposition of $End(W)$ in \eqref{od}.
Hence for $t\, \in\, \mathbb R$, all four $tv$, $tv^*$,
$tw$ and $tw^*$ are $1$--forms with values in
$F_0$. On the other hand, $\theta$ and $\theta^*$ are
$1$--forms with values in $\text{ad}(E^0_G)$.

Let $\nabla'$ be the connection on the vector bundle $W$ associated
to $E^0_G$ induced by the connection $\nabla$ on $E^0_G$. So the
curvature of $\nabla'$ coincides with the curvature of $\nabla$,
in particular, $K(\nabla')$ is a $2$--form with values in 
$\text{ad}(E^0_G)$. We note that
$$
K(\nabla'+tv - tv^*)\, =\, K(\nabla') +t\nabla'(v-v^*) +t^2C\, ,
$$
where $C$ is independent of $t$. Since $\nabla'$ is induced
by a connection $E^0_G$, and the decomposition in \eqref{od2}
is preserved by the action of $G$, the connection
$\nabla'$ preserves the decomposition in \eqref{od}. Hence
$\nabla'(v-v^*)$ is a $2$--form with values in $F_0$.

Using these and the fact that the decompositions in \eqref{o1}
are orthogonal, we have
$$
(\frac{d}{dt}
\Vert K(\nabla'+tv - tv^*)+[\theta+tw\, ,
\theta^*+tw^*]\Vert^2)\Big\vert_{t=0}\,=\, 0\,.
$$
This completes the proof of the lemma.
\end{proof}

Let
\begin{equation}\label{P}
\Psi_W\, :\, {\mathcal B}_W\, \longrightarrow\,
T^{\mathbb R} {\mathcal B}_W
\end{equation}
be the gradient vector field on
${\mathcal B}_W$ for the function $\text{YMH}_W$. From
Lemma \ref{lem2} and \eqref{i} we have the following
corollary:

\begin{corollary}\label{cor1}
The restriction of the vector field $\Psi_W$ to $\phi({\mathcal
B}_0)$ (see \eqref{e9}) lies in the image of the differential
$d\phi$ in \eqref{e10}. Furthermore, this restriction
coincides with the gradient vector field for the function 
${\rm YMH}_G$.
\end{corollary}

\section{Closedness of the embedding}\label{sec:closed}

For a complex vector space $V'$, let $P(V')$ denote the
projective space of lines in $V'$. Any linear action on
$V'$ induces an action on $P(V')$.

Consider the closed
subgroup $\rho(G)\, \subset\, \text{GL}(V)$ in \eqref{e2}.
A theorem of C. Chevalley (see
\cite[p. 80]{Hu}) says that there is a finite dimensional left
$\text{GL}(V)$--module $V_1$ and a line
\begin{equation}\label{l}
\ell\, \subset\, V_1
\end{equation}
such that $\rho(G)$ is exactly the isotropy subgroup,
for the action of $\text{GL}(V)$ on $P(V_1)$,
of the point in $P(V_1)$ representing the line $\ell$.

Let $E_{\text{GL}(V)}\, :=\, E^0_G(\text{GL}(V))
\, =\, E^0_G\times_G \text{GL}(V)\, \longrightarrow\, X$
be the principal $\text{GL}(V)$--bundle obtained by extending
the structure group of $E^0_G$ (see \eqref{e4}) by the homomorphism
$\rho$ in \eqref{e2}. Therefore, the vector bundle
$E_{\text{GL}(V)}(V)$, associated to $E_{\text{GL}(V)}$ by the
standard action of $\text{GL}(V)$ on $V$, is identified with
the vector bundle $W$ in \eqref{W}. Let
\begin{equation}\label{v1}
{\mathcal V}_1 \, :=\, E_{\text{GL}(V)}(V_1)\, \longrightarrow\, X
\end{equation}
be the vector bundle associated to $E_{\text{GL}(V)}$ for the
above $\text{GL}(V)$--module $V_1$. Since
$$
{\mathcal V}_1\, =\, E^0_G(V_1)\, ,
$$
and the action of $G$ on $V_1$ preserves the line $\ell$
in \eqref{l}, the line $\ell$ defines a $C^\infty$ line subbundle
\begin{equation}\label{l0}
L_0\, \subset\, {\mathcal V}_1\, .
\end{equation}

Take any holomorphic structure $\overline{\partial}_W
\, \in\,{\mathcal A}(W)$ on the
vector bundle $W$ (see \eqref{eb1}). The holomorphic structure 
$\overline{\partial}_W$ on $W$ defines a holomorphic structure
on the principal $\text{GL}(V)$--bundle $E_{\text{GL}(V)}$
corresponding to $W$. Hence $\overline{\partial}_W$ defines a
holomorphic structure on
the associated vector bundle ${\mathcal V}_1$
in \eqref{v1}. This holomorphic structure on ${\mathcal V}_1$
will be denoted by $\overline{\partial}'_1$.

Since $\rho(G)$ is the isotropy subgroup of the point in
$P(V_1)$ representing the line $\ell$ for the action of
$\text{GL}(V)$ on $P(V_1)$, we conclude that
$\overline{\partial}_W$ lies in $\delta({\mathcal A}_0)$
(see \eqref{eb1}) if and only if the line subbundle $L_0$
in \eqref{l0} is holomorphic with respect to the holomorphic
structure $\overline{\partial}'_1$ on ${\mathcal V}_1$.

Therefore, we have the following lemma:

\begin{lemma}\label{lem3}
The embedding $\delta$ in \eqref{eb1} is closed.
\end{lemma}

The action of $\text{GL}(V)$ on $V_1$ gives a homomorphism
$$
\text{End}_{\mathbb C}(V)\, \longrightarrow\, \text{End}_{\mathbb 
C}(V_1)
$$
of Lie algebras. This homomorphism in turn gives a
homomorphism of vector bundles
\begin{equation}\label{f1}
End(W)\, \longrightarrow\, End({\mathcal V}_1)\, ,
\end{equation}
where ${\mathcal V}_1$ is the vector bundle in \eqref{v1}.

Take any $\theta\, \in\, A^{1,0}(End(W))$. Let
$\theta'\, \in\, A^{1,0}(End({\mathcal V}_1))$ be the section
constructed from $\theta$ using the homomorphism in \eqref{f1}.
Since $\rho(G)$ is the isotropy subgroup of the point in
$P(V_1)$ representing the line $\ell$ for the action of
$\text{GL}(V)$ on $P(V_1)$, we conclude the following:
The section $\theta$ lies in the image of the natural homomorphism
$$
A^{1,0}(\text{ad}(E^0_G))\, \longrightarrow\,
A^{1,0}(End(W))
$$
if and only if $\theta'(L_0)\, \in\, A^{1,0}(L_0)$, where
$L_0$ is the subbundle in \eqref{l0}.

Therefore, using Lemma \ref{lem3}, we have following proposition:

\begin{proposition}\label{prop1}
The embedding $\phi$ in \eqref{e9} is closed.
\end{proposition}

Let
\begin{equation}\label{cs}
{\mathcal S}(E^0_G)\, \subset\, {\mathcal B}_0
\end{equation}
be the subset consisting of all pairs that are Higgs
$G$--bundles. So a pair
$$
(\overline{\partial}_{E_G^0}\, ,\theta)\,
\in\, {\mathcal A}_0\times A^{1,0}(\text{ad}(E^0_G))\,=\, 
{\mathcal B}_0
$$
lies in ${\mathcal S}(E^0_G)$ if and only if the
section $\theta$ is holomorphic with respect to the holomorphic 
structure on $\text{ad}(E^0_G)\otimes K_X$ defined by 
$\overline{\partial}_{E_G^0}$.

Consider the gradient flow on ${\mathcal B}_0$ for the function
$\text{YMH}_G$ defined in \eqref{e11}. The following lemma shows
that this flow preserves the subset ${\mathcal S}(E^0_G)$ defined
in \eqref{cs}.

\begin{lemma}\label{lem4}
Take any $z\, :=\, (\overline{\partial}_{E_G^0}\, ,\theta)\,
\in\, {\mathcal S}(E^0_G)$. Let
$$
\gamma_z\, :\, {\mathbb R}\, \longrightarrow\,{\mathcal B}_0
$$
be the integral curve for the gradient flow on ${\mathcal B}_0$
for the function ${\rm YMH}_G$ such that $\gamma_z(0)\,=\,z$.
Then
$$
\gamma_z(t)\, \in\, {\mathcal S}(E^0_G)
$$
for all $t\, \in\, {\mathbb R}$.
\end{lemma}

\begin{proof}
Consider ${\mathcal B}_W$ defined in \eqref{e8}. Let
$$
{\mathcal S}(W)\, \subset\, {\mathcal B}_W
$$
be the subset consisting
of all pairs $(\overline{\partial}'\, ,\theta)\,\in
\,{\mathcal A}(W)\times A^{1,0}(End(W))$ such that $\theta$ is
holomorphic with respect to the holomorphic structure given by
$\overline{\partial}'$. For the map $\phi$ in \eqref{e9},
$$
\phi({\mathcal S}(E^0_G))\, \subset\, {\mathcal S}(W)
$$
(see Remark \ref{rem1}).

In view of Corollary \ref{cor1}, to prove the lemma it suffices
to show that the vector field $\Psi_W$ (defined in \eqref{P})
preserves the subset ${\mathcal S}(W)$. But this is proved in
\cite{Wi}; from \cite[Lemma 3.10]{Wi} and the proof
of Proposition 3.2 in \cite[pp. 295--297]{Wi} it follows that the
flow $\Psi_W$ is generated by the action of the complex gauge
group, hence ${\mathcal S}(W)$ is preserved by the flow.
This completes the proof of the lemma.
\end{proof}

\begin{theorem}\label{lw}
The integral curve $\gamma_z$ for the gradient flow of 
${\rm YMH}_G$ on $\mathcal{B}_0$ with initial condition $z \, :=\, 
(\overline{\partial}_{E_G^0}\, ,\theta)\, \in\, {\mathcal S}(E^0_G)$ 
converges to a limit in ${\mathcal S}(E^0_G)$.
\end{theorem}

\begin{proof}
Theorem 1.1 in \cite{Wi} shows that the gradient flow of ${\rm YMH}_W$ 
on the space $\mathcal{B}_W$ with initial conditions in $\mathcal{S}(W)$ 
converges to a limit in $\mathcal{S}(W)$. Moreover, 
Corollary \ref{cor1} and Lemma \ref{lem4} together with the uniqueness 
of the flow from Proposition 3.2 in \cite{Wi} give the
following: when 
the initial conditions are in $\phi \left( \mathcal{S}(E^0_G) \right)$, 
then the flow preserves the space $\phi \left( \mathcal{S}(E^0_G) 
\right)$. Combining these two facts, we see that because the embedding 
$\phi$ is closed by Proposition \ref{prop1}, the limit of the flow 
lies in $\phi \left( \mathcal{S}(E^0_G) \right)$. Since $\phi(\gamma_z)$ 
coincides with the gradient flow of ${\rm YMH}_W$ with initial 
conditions in $\phi \left( \mathcal{S}(E^0_G) \right)$ by Corollary 
\ref{cor1}, we conclude that $\lim_{t \rightarrow \infty} 
\gamma_z(t)$ exists, and it is in $\mathcal{S}(E^0_G)$.
\end{proof}

\section{Reduction of structure group}\label{sec:reduction}

As before, $G$ is a connected reductive linear algebraic group
defined over $\mathbb C$.

See \cite{BG}, \cite{BiSc} for the definitions
of semistable, stable and polystable Higgs $G$--bundles.

\begin{lemma}\label{lem5}
Let $(E_G\, ,\theta)$ be a semistable Higgs $G$--bundle on $X$.
The Higgs vector bundle $({\rm ad}(E_G)\, ,\varphi)$ is
semistable, where $\varphi$ is the Higgs field on ${\rm ad}(E_G)$
defined by $\theta$ using the Lie algebra structure of the
fibers of ${\rm ad}(E_G)$.
\end{lemma}

\begin{proof}
This follows from \cite[p. 37, Lemma 3.6]{BG}, but some
explanations are necessary.

Let $Z(G)\, \subset\,G$ be the connected component of the center
of $G$ containing the identity element. Define
$$
G'\, :=\, G/Z(G)\, .
$$
Let $(E_{G'}\, ,\theta')$ be the Higgs $G'$--bundle
over $X$ obtained by
extending the structure group of $(E_G\, ,\theta)$ using the
quotient map $G\, \longrightarrow\, G'$. Since $(E_G\, ,\theta)$
is semistable, it follows immediately that
the Higgs $G'$--bundle $(E_{G'}\, ,\theta')$ is semistable.
Let $\varphi'$ be the Higgs field on the adjoint vector
bundle $\text{ad}(E_{G'})$ induced by $\theta'$. Since
$(E_{G'}\, ,\theta')$ is semistable, and the group $G'$ does
not have any nontrivial character, the Higgs vector
bundle $(\text{ad}(E_{G'})\, ,\varphi')$ is semistable
\cite[p. 37, Lemma 3.6]{BG} (see also
\cite[p. 26, Proposition 2.4]{BG}). We have
$$
({\rm ad}(E_G)\, ,\varphi)\,=\, (\text{ad}(E_{G'})\, ,\varphi')\oplus
(X\times z({\mathfrak g}), 0)\, ,
$$
where $z({\mathfrak g})$ is the Lie algebra of
$Z(G)$, and $X\times z({\mathfrak g})$ is the trivial vector bundle
over $X$ with
fiber $z({\mathfrak g})$. Hence $({\rm ad}(E_G)\, ,\varphi)$ is
semistable; note that $\text{degree}({\rm ad}(E_G))\,=\, 0\,=\,
\text{degree}({\rm ad}(E_{G'}))$. This completes the proof of
the lemma.
\end{proof}

Let $(E_G\, ,\theta)$ be a Higgs $G$--bundle on $X$ and
$H\,\subset\, G$ a closed algebraic subgroup. A 
\textit{reduction} of structure group of the Higgs $G$--bundle
$(E_G\, ,\theta)$ to $H$ is a holomorphic reduction of
structure group $E_H \,\subset\, E_G$ to $H$ over $X$ such that
$\theta$ lies in the image of the
homomorphism $H^0(X,\, {\rm ad}(E_H)\otimes K_X)
\, \longrightarrow\, H^0(X,\, {\rm ad}(E_G)\otimes K_X)$.

Given a Higgs $G$--bundle $(E_G\, ,\theta)$, there is a canonical
Harder--Narasimhan reduction of structure group of
$(E_G\, ,\theta)$ to a parabolic subgroup $P$ of $G$ \cite{DP}
(the method in \cite{DP} is based on \cite{BH}). If
$(E_G\, ,\theta)$ is semistable, then $P\,=\, G$.

We recall the definition of the Harder--Narasimhan reduction
of a Higgs $G$--bundle.

Let $(E_G\, ,\theta)$ be a Higgs $G$--bundle on $X$. Then there
is a parabolic subgroup $P\, \subset\, G$ and a reduction
of structure group $E_P$ of $(E_G\, ,\theta)$ to $P$ such that
\begin{enumerate}
\item the principal $L(P)$--bundle $E_{L(P)}\, :=\, E_P\times_P
L(P)\, \longrightarrow\, X$, where $L(P)$ is the Levi quotient 
of $P$, is semistable, and

\item for any nontrivial character $\chi$ of $P$ which is
a nonnegative linear combination of simple roots (with
respect to some Borel subgroup contained in $P$) and is
trivial on the center of $G$, the associated line bundle
$E_P(\chi)\, \longrightarrow\, X$ is of positive degree.
\end{enumerate}
The above pair $(P\, ,E_P)$ is unique in the following sense:
for any other pair $(P_1\, , E_{P_1})$ satisfying
the above two conditions, there is some $g\, \in\, G$ such that
\begin{itemize}
\item $P_1\, =\, g^{-1}Pg$, and

\item $E_{P_1}\, =\, E_Pg$.
\end{itemize}
(See \cite{DP}, \cite{BH}.)

A semistable vector bundle $E\, \longrightarrow\, X$
admits a filtration of subbundles
\begin{equation}\label{s1}
0 \, =\, E_0\, \subset\, E_1 \, \subset\, E_2\, \subset\,
\cdots\, \subset\, E_{n-1} \, \subset\, E_n\, =\, E
\end{equation}
such that $E_i/E_{i-1}$, $1\,\leq\, i\, \leq\, n$, is the maximal
polystable subbundle of $E/E_{i-1}$ with
$$
\frac{\text{degree}(E_i/E_{i-1})}{\text{rank}(E_i/E_{i-1})}\,=\,
\frac{\text{degree}(E)}{\text{rank}(E)}
$$
(see \cite[p. 23, Lemma 1.5.5]{HL});
this filtration is called the \textit{socle filtration}.
In \cite{AB}, this was generalized to semistable principal
$G$--bundles (see \cite[p. 218, Proposition 2.12]{AB}). In Theorem
\ref{thm1} proved below, this is further generalized to
semistable Higgs $G$--bundles.

We will define admissible reductions of a Higgs $G$--bundle.
See \cite[pp. 3998--3999]{BS} for the definition of an admissible 
reduction of structure group of a principal $G$--bundle.

\begin{definition}\label{defa}
{\rm An {\em admissible} reduction of structure group of
a Higgs $G$--bundle $(E_G\, ,\theta)$ to a parabolic
subgroup $P\, \subset\, G$ is a
reduction of structure group $E_P$ of $(E_G\, ,\theta)$ to
$P$ such that $E_P\, \subset\, E_G$ is an
admissible reduction of $E_G$.}
\end{definition}

Let $(E_G\, ,\theta)$ be a Higgs $G$--bundle on $X$. Let
$E'_P\, \subset\, E_G$ be a reduction of structure group of
$(E_G\, ,\theta)$ to a parabolic subgroup $P$ of $G$. So
$\theta$ is a section of $\text{ad}(E'_P)\otimes K_X$.
Let $L(P)$ be the Levi quotient of $P$. Let
$E'_P(L(P))$ be the principal
$L(P)$--bundle over $X$ obtained by extending the structure group
of $E'_P$ using the quotient map $P\, \longrightarrow\, L(P)$. The 
quotient homomorphism $\text{Lie}(P)\, \longrightarrow\, 
\text{Lie}(L(P))$ induces a homomorphism of adjoint bundles
$$
\text{ad}(E'_P)\, \longrightarrow\, \text{ad}(E'_P(L(P)))\, .
$$
Using this homomorphism of vector bundles, the section $\theta$
of $\text{ad}(E'_P)\otimes K_X$ gives a holomorphic section of
$\text{ad}(E'_P(L(P)))\otimes K_X$. In other words, $\theta$ gives
a Higgs field on $E'_P(L(P))$. This Higgs field on $E'_P(L(P))$
will be denoted by $\theta'$.

Let $(E_G\, ,\theta)$ be a semistable Higgs
$G$--bundle on $X$ which is not polystable.
Let $Q\, \subsetneq\, G$ be a proper parabolic subgroup which is
maximal among all the proper parabolic
subgroups $P$ such that $(E_G\, ,\theta)$ has an
admissible reduction of structure group
$E'_P\, \subset\, E_G$
(see Definition \ref{defa}) for which the associated Higgs
$L(P)$--bundle $(E'_P(L(P))\, ,\theta')$ defined above
is polystable.

\begin{definition}\label{g-s}
{\rm An admissible reduction of structure group of
$(E_G\, ,\theta)$ to $Q$
$$
E_Q\, \subset\, E_G
$$
will be called a {\em socle reduction} if the associated
Higgs $L(Q)$--bundle $(E_Q(L(Q))\, ,\theta')$ is polystable,
where $L(Q)$ is the Levi quotient of $Q$.}
\end{definition}

\begin{theorem}\label{thm1}
Let $(E_G\, ,\theta)$ be a semistable Higgs
$G$--bundle on $X$ which is not polystable. Then 
$(E_G\, ,\theta)$ admits a socle reduction. If
$(Q\, , E_Q)$ and $(Q_1\, , E_{Q_1})$ are two
socle reductions of $(E_G\, ,\theta)$, 
then there is some $g\, \in\, G$ such that
$Q_1\, =\, g^{-1}Qg$, and $E_{Q_1}\, =\, E_Qg$.
\end{theorem}

\begin{proof}
First note that the construction of the socle filtration of a
semistable vector bundle extends to semistable Higgs bundles;
indeed, the proof in \cite[p. 23, Lemma 1.5.5]{HL} goes through in this
case also. Therefore, if $(E\, ,\theta)$ is a semistable Higgs
vector bundle on $X$ which is not polystable, there is filtration
of subbundles
\begin{equation}\label{f}
0\, \subset\, E_1 \, \subset\, E_2\, \subset\,\cdots
\, \subset\, E_{n-1} \, \subset\, E_n\, =\, E
\end{equation}
such that
\begin{equation}\label{f2}
\theta(E_i)\, \subset\, E_i\otimes K_X
\end{equation}
for all
$i\, \in\, [1\, ,n]$, and $(E_i/E_{i-1}\, ,\theta'_i)$
is the unique maximal polystable Higgs subbundle of the
Higgs bundle $(E/E_{i-1}\, ,\theta''_i)$ such that
$$
\frac{\text{degree}(E_i/E_{i-1})}{\text{rank}(E_i/E_{i-1})}\,=\,
\frac{\text{degree}(E)}{\text{rank}(E)}\, ,
$$
where $\theta'_i$ and
$\theta''_i$ are the Higgs fields on $E_i/E_{i-1}$ and $E/E_{i-1}$
respectively induced by $\theta$ (the condition in \eqref{f2}
ensures that $\theta$ induces Higgs fields on both
$E_i/E_{i-1}$ and $E/E_{i-1}$).

Let $\text{ad}(E_G)\, \longrightarrow\, X$ be the adjoint bundle
of $E_G$. Let $\varphi$ be the Higgs field on $\text{ad}(E_G)$
defined by $\theta$. From Lemma
\ref{lem5} we know that the Higgs vector
bundle $(\text{ad}(E_G)\, ,\varphi)$ is semistable.
We note if $(\text{ad}(E_G)\, ,\varphi)$ is polystable, then
$(E_G\, ,\theta)$ is polystable. Since $(E_G\, ,\theta)$
is not polystable, we conclude that the Higgs
vector bundle $(\text{ad}(E_G)\, ,\varphi)$ is not polystable.
Let
\begin{equation}\label{f4}
0\, =\, E_0 \, 
\subset\, E_1 \, \subset\, E_2\, \subset\,\cdots
\, \subset\, E_{m-1} \, \subset\, E_m\, =\, \text{ad}(E_G)
\end{equation}
be the socle filtration for $(\text{ad}(E_G)\, ,\varphi)$
(see \eqref{f}).

Fix a $G$--invariant nondegenerate symmetric bilinear form $B_0$ on
the Lie algebra $\mathfrak g$ of $G$; such a form exists because $G$
is reductive. This form $B_0$ defines nondegenerate symmetric bilinear
forms on the fibers of $\text{ad}(E_G)$. So we get an isomorphism
\begin{equation}\label{is}
\text{ad}(E_G)\, \stackrel{\sim}{\longrightarrow}\,
\text{ad}(E_G)^*\, .
\end{equation}
Let $\varphi^*$ be the dual Higgs field on $\text{ad}(E_G)^*$ defined
by $\varphi$. The isomorphism in \eqref{is} clearly takes $\varphi$
to $\varphi^*$. In particular, $(\text{ad}(E_G)\, ,\varphi)$ is
self--dual.

{}From the uniqueness of the socle filtration it follows
that the filtration in \eqref{f4} is self--dual. Also,
the integer $m$ in \eqref{f4} is odd. The tensor product of two 
semistable Higgs bundles on $X$ is again semistable
\cite[p. 38, Corollary 3.8]{Si1}. Using these observations
it follows that
\begin{itemize}
\item the subbundle $E_{\frac{m+1}{2}}$ is closed under the Lie
bracket operation on the fibers of $\text{ad}(E_G)$,

\item the fibers of $E_{\frac{m-1}{2}}$ are ideals in the
fibers of $E_{\frac{m+1}{2}}$, and are nilpotent,

\item the fibers of the quotient $E_{\frac{m+1}{2}}/
E_{\frac{m-1}{2}}$ are reductive, and

\item the Higgs field $\theta$ is a section of
$E_{\frac{m+1}{2}}\otimes K_X$.
\end{itemize}
(see \cite[p. 218, Proposition 2.12]{AB}). It should be clarified
that to prove the above statements we need the following: for any two
polystable Higgs vector bundles $(W_1\, , \varphi_1)$ and 
$(W_2\, , \varphi_2)$ over $X$, the tensor product $(W_1\otimes W_2\, ,
\varphi_1\otimes\text{Id}_{W_2}+ \text{Id}_{W_1}\otimes\varphi_2)$
is also a polystable Higgs vector bundle. To prove that
$(W_1\otimes W_2\, ,
\varphi_1\otimes\text{Id}_{W_2}+ \text{Id}_{W_1}\otimes\varphi_2)$
is polystable, let $\nabla_1$ and
$\nabla_2$ be the Hermitian--Yang--Mills connections on
$(W_1\, , \varphi_1)$ and $(W_2\, , \varphi_2)$ respectively
(see \cite[p. 19, Theorem 1(2)]{Si1}).
Then the induced connection $\nabla_1\otimes\text{Id}_{W_2}+ 
\text{Id}_{W_1}\otimes\nabla_2$ on $W_1\otimes W_2$ is
a Hermitian--Yang--Mills connection for
$(W_1\otimes W_2\, ,
\varphi_1\otimes\text{Id}_{W_2}+ \text{Id}_{W_1}\otimes\varphi_2)$.
Hence $(W_1\otimes W_2\, ,
\varphi_1\otimes\text{Id}_{W_2}+ \text{Id}_{W_1}\otimes\varphi_2)$
is polystable \cite[p. 19, Theorem 1(2)]{Si1}.

{}From the above statements it follows that
$E_{\frac{m+1}{2}}$ is a Lie algebra subbundle of
the Lie algebra bundle $\text{ad}(E_G)$ such that the
fibers of $E_{\frac{m+1}{2}}$ are parabolic subalgebras. 

The fibers of $\text{ad}(E_G)$ are identified with the Lie algebra
$\mathfrak g$ up to an inner automorphism. More precisely, for
any point $x\, \in\, X$, and any point $z$ in the fiber $(E_G)_x$
of $E_G$, we have an isomorphism
\begin{equation}\label{il}
\sigma_z \, :\, {\mathfrak g}\, \longrightarrow\, \text{ad}(E_G)_x
\end{equation}
that sends any $v\, \in\, {\mathfrak g}$ to the
image of $(z\, ,v)$ in $\text{ad}(E_G)_x$ (recall that
$\text{ad}(E_G)$ is a quotient of $E_G\times{\mathfrak g}$).
For any $g\, \in\, G$, the isomorphisms $\sigma_z$ and
$\sigma_{zg}$ differ by the inner automorphism
$\text{Ad}(g)$ of ${\mathfrak g}$.
Let $Q\, \subset\, G$ be a parabolic subgroup in the conjugacy
class of parabolic subgroups whose Lie algebras are
identified with the fibers of $E_{\frac{m+1}{2}}$
by some isomorphism constructed
in \eqref{il}. The normalizer of any parabolic
subgroup $P\, \subset\, G$ coincides with $P$. In particular,
the normalizer of $Q\, \subset\, G$ is $Q$ itself. Hence
the subalgebra
bundle $E_{\frac{m+1}{2}}\, \subset\, \text{ad}(E_G)$ gives
a holomorphic reduction of structure group $E_Q\, \subset\, E_G$
such that the subbundle $\text{ad}(E_Q)\, \subset\, \text{ad}(E_G)$
coincides with $E_{\frac{m+1}{2}}$. For any point $x\, \in\, X$,
the fiber $(E_Q)_x\, \subset\, (E_G)_x$ consists of all points
$z\, \in\, (E_G)_x$
such that the isomorphism $\sigma_z$ in \eqref{il} takes
$\text{Lie}(Q)$ to $(E_{\frac{m+1}{2}})_x$.

Since $\text{ad}(E_Q)\,=\, E_{\frac{m+1}{2}}$, and the Higgs field
$\theta$ is a section of $E_{\frac{m+1}{2}}\otimes K_X$, we conclude
that $E_Q$ is a reduction of structure group of the Higgs $G$--bundle
$(E_G\, ,\theta)$. It is straight--forward to check that $E_Q$
is a socle reduction of $(E_G\, ,\theta)$.

Given any socle reduction $E_{Q'}$ of $(E_G\, ,\theta)$, it can be
shown that the adjoint bundle $\text{ad}(E_{Q'})$ coincides
with the subbundle $E_{\frac{m+1}{2}}$ in \eqref{f4}. From this
the uniqueness statement in the theorem follows.
This completes the proof of the theorem.
\end{proof}

For a polystable Higgs $G$--bundle $(E_G\, ,\theta)$ the 
\textit{socle reduction} is defined to be $E_G$ itself.

Given a Higgs $G$--bundle, combining the Harder--Narasimhan reduction
with the socle reduction we get a new Higgs $G$--bundle which will
be described below.

Let $(E_G\, ,\theta)$ be a Higgs $G$--bundle. Let
$(E_P\, ,\theta_P)$ be the Harder--Narasimhan reduction of
$(E_G\, ,\theta)$. If $(E_G\, ,\theta)$ is semistable, then
$P\, =\, G$, and $(E_P\, ,\theta_P)\,=\, (E_G\, ,\theta)$.

Let $L(P)$ be the Levi quotient of $P$. Let
\begin{equation}\label{elp}
(E_{L(P)}\, ,\theta_{L(P)})
\end{equation}
be the Higgs $L(P)$--bundle obtained by
extending the structure group of the above Higgs $P$--bundle $(E_P\, 
,\theta_P)$ using the quotient map $P\, \longrightarrow\, L(P)$. From
the definition of a Harder--Narasimhan reduction we know that
the Higgs $L(P)$--bundle $(E_{L(P)}\, ,\theta_{L(P)})$ is semistable.
Therefore, $(E_{L(P)}\, ,\theta_{L(P)})$ has a unique socle reduction
by Theorem \ref{thm1}. Let
$$
E_H\, \subset\, E_{L(P)}
$$
be the socle reduction of $(E_{L(P)}\, ,\theta_{L(P)})$.
So $H$ is a Levi subgroup of a parabolic subgroup
of $L(P)$; the
Higgs field on $E_H$ induced by $\theta_{L(P)}$ will
be denoted by $\theta_H$ (see Definition \ref{g-s}).

The Levi quotient $L(P)$ is identified with all the Levi factors
of $P$, and $H$ is a subgroup of $L(P)$. Therefore, $H$ becomes
a subgroup of $G$ after fixing a Levi factor of $P$. Let
\begin{equation}\label{ep}
(E'_G\, ,\theta')
\end{equation}
be the Higgs $G$--bundle obtained
by extending the structure group of the Higgs $H$--bundle $(E_H\, 
,\theta_H)$ using the inclusion of $H$ in $G$.

Take any
\begin{equation}\label{z}
z \, :=\, (\overline{\partial}_{E_G^0}\, ,\theta)\, \in\, {\mathcal 
S}(E^0_G)
\end{equation}
(see \eqref{cs}). Let $(E_G\, ,\theta)$ be the
Higgs $G$--bundle defined by the holomorphic structure
$\overline{\partial}_{E_G^0}$ on $E^0_G$ together with the section 
$\theta$
in \eqref{z}. Let $(E'_G\, ,\theta')$ be the new Higgs $G$--bundle
constructed in \eqref{ep} from $(E_G\, ,\theta)$.

\begin{lemma}\label{lem6}
Let $\gamma_z$ be the integral curve for the gradient flow of
${\rm YMH}_G$ on $\mathcal{B}_0$ with initial condition $z$
(see \eqref{z}). Let
$$
(\overline{\partial}_1\, ,\theta_1)\, =\, \lim_{t\to \infty}
\gamma_z(t)\, \in\, {\mathcal S}(E^0_G)
$$
be the limit in Theorem \ref{lw}. Then the Higgs $G$--bundle
defined by $(\overline{\partial}_1\, ,\theta_1)$ is
holomorphically isomorphic to the Higgs $G$--bundle
$(E'_G\, ,\theta')$ constructed above.
\end{lemma}

\begin{proof}
For Higgs vector bundles this was proved in \cite{Wi}
(see \cite[p. 325, Theorem 5.3]{Wi}). Let $(E_G\, ,\theta)$
be the Higgs $G$--bundle defined by $z$ in \eqref{z}.
Let $({\rm ad}(E_G)\, ,\varphi)$ be the corresponding 
Higgs vector bundle defined by the Higgs field on the
adjoint vector bundle ${\rm ad}(E_G)$ induced by $\theta$.

Recall that the Harder--Narasimhan reduction of the
Higgs $G$--bundle $(E_G\, ,
\theta)$ is constructed using the Harder--Narasimhan filtration of
the Higgs vector bundle
$({\rm ad}(E_G)\, ,\varphi)$. Let $(E_{L(P)}\, ,\theta_{L(P)})$
be the semistable principal Higgs bundle constructed as in
\eqref{elp} from the Harder--Narasimhan reduction of $(E_G\, ,
\theta)$. Recall that the socle reduction of
a semistable Higgs $L(P)$--bundle $(E_{L(P)}\, ,\theta_{L(P)})$
is constructed using the socle filtration of the
adjoint vector bundle $\text{ad}(E_{L(P)})$
equipped with the Higgs field induced by $\theta_{L(P)}$. From
these constructions it can be deduced that the
Harder--Narasimhan--socle filtration of the Higgs vector bundle
$({\rm ad}(E_G)\, ,\varphi)$ is compatible with the filtration
of ${\rm ad}(E_G)$ obtained from $(E_{L(P)}\, ,\theta_{L(P)})$.
Using this, the lemma follows.
\end{proof}


\end{document}